\documentclass[preprint,review,12pt]{elsarticle}
\usepackage{amsthm}
\usepackage{tikz}         % for drawing graphics
\usepackage[numbers]{natbib}
\setcitestyle{square,numbers,sort&compress}
\usepackage[utf8]{inputenc}
\usepackage[T1]{fontenc}
\usepackage{hyperref}       % hyperlinks
\usepackage{url}            % simple URL typesetting
\usepackage{booktabs}       % professional-quality tables
\usepackage{tabularx}       % tables with automatic width adjustment and line wrapping
\usepackage{amsmath}        % mathematical typesetting
\usepackage{amsfonts}       % blackboard math symbols
\usepackage{nicefrac}       % compact symbols for 1/2, etc.
\usepackage{microtype}      % microtypography
\usepackage{lipsum}		% Can be removed after putting your text content
\usepackage{graphicx}
\usepackage{lineno}
\usepackage{doi}
\usepackage{bm}
\usepackage{float}
\setcitestyle{square}
\usepackage[a4paper,left=2.0cm,right=2.0cm,top=2.4cm,bottom=2.4cm]{geometry}

\newtheorem{theorem}{Theorem}

\newtheorem{assumption}{Assumption}

\begin{document}

\begin{frontmatter}
\title{Multiresolution Sobolev Trained Fourier Neural Operator for Cross-Resolution Surrogate Modelling of Parametric Two-Phase Darcy Flows in Porous Media}

\author[1,2]{Zhao Zhang\corref{cor1}}
\ead{zhaozhang@sdu.edu.cn}
\author[1,2]{Chengjin Guo}
\author[1,2]{Zhenglong Chen}
\author[4,5]{Kai Zhang}
\author[5]{Piyang Liu}
\author[4]{Xia Yan\corref{cor1}}
\ead{jsyanxia1989@163.com}
\address[1]{Research Centre for Mathematics and Interdisciplinary Sciences, Shandong University, Qingdao, Shandong Province, 266237, China}
\address[2]{Frontiers Science Center for Nonlinear Expectations, Minister of Education, Shandong University, Qingdao, Shandong Province, 266237, China}
\address[4]{School of Petroleum Engineering, China University of Petroleum (East China), Qingdao, Shandong Province, 266580, China}
\address[5]{School of Civil Engineering, Qingdao University of Technology, Qingdao, Shandong Province, 266520, China}
\cortext[cor1]{Corresponding author}

\begin{abstract}
Surrogate modelling is important to enhance the computational efficiency of uncertainty quantification for physical processes involving random parameters. The advantage of Fourier neural operator (FNO) compared to conventional surrogates is its cross-resolution prediction capability. In the current study, a cross-resolution surrogate based on FNO is built for two-phase transient Darcy flows in porous media with random parameters. A multiresolution Sobolev learning method is proposed for FNO to enhance the cross-resolution prediction accuracy by incorporating first and second-order gradients in the loss function using training samples of different resolutions. The gradients are approximated by resolution-related finite difference schemes compatible with the grid for numerical simulation. Experimental results demonstrate the effectiveness of the new method in training FNO for building cross-resolution surrogates.

\end{abstract}

\end{frontmatter}

\newcommand{\h}{\hat}
%\linenumbers
\section{Introduction}
The objective in many science and engineering problems is not to solve a single governing partial differential equation (PDE), but to learn a solution operator that maps from the uncertain model parameters, such as physical properties, boundary conditions, or source terms, to the solution fields \citep{azizzadenesheli2024neural,brunton2024promising}. The learned operator can serve as a surrogate to enhance the computational efficiency of forward and inverse uncertainty quantification \citep{smith2013uncertainty, zhang2024ensemble}.
In recent years, neural operators, such as DeepONet\citep{Lu2021-oy} and Fourier Neural Operators (FNO)\citep{li2020fourier} to name a few, have been proposed for operator learning \citep{lu2022comprehensive}. Compared to conventional surrogate modelling methods, a major advantage of neural operators is the cross-resolution generalization capability \citep{li2024physics}.

Neural operators can be data-driven or physics-informed \citep{karniadakis_physics-informed_2021,liu2025architectures}. Data-driven learning requires sufficient training data while physics-informed learning can be conducted with limited or no labelled data in some cases \citep{zhang2022physics,zhang2023physics,yan2024physics, yan2025finite}. However, it is still difficult for physics-informed learning to converge even after extensive hyperparameter tuning in the operator learning of two-phase compressible Darcy flow cases involving time-dependent nonlinear PDEs and heterogeneous parameter fields. 

In the current study, we build a data-driven cross-resolution surrogate based on FNO for two-phase Darcy flows in porous media involving random parameters. A multiresolution Sobolev learning scheme is proposed for FNO on training samples of different resolutions.  Using samples of lower resolutions that are cheaper to be generated, the trained surrogate can be used to predict flow fields of higher resolution. We show that using multiresolution training samples can effectively enhance the cross-resolution prediction accuracy of FNO-based surrogates. A multiresolution Sobolev training algorithm is established to further enhance the cross-resolution prediction accuracy by incorporating gradients into the loss function. The gradients are approximated by resolution-related finite difference schemes compatible with the grid.
Conventional Sobolev learning often involves only first-order gradients in the loss function \citep{tsay2021sobolev,yuan2022sobolev,cho2025sobolev,oh2025sobolev}. In this paper, both first and second-order gradients are accounted for. Gradient constraints can be regarded as regularizers to help restrict the learning process to reach better lower minima. 

The paper is organized as follows. In section 2, the problem background of parametric Darcy flow and operator learning is introduced. In section 3, the numerical methods for data generation, FNO structure and loss function are presented. In section 4, the numerical experiments and results are discussed. In section 5, the conclusions are presented.

\section{Problem Setting}
\subsection{Parametric Two-phase Compressible Darcy Flow}
In the current study, operator learning is applied to Darcy flows in porous media for training cross-resolution surrogates with limited labelled data. 
The nonlinear PDEs governing immiscible two-phase compressible Darcy flows in porous media are \citep{chen2006computational}
\begin{eqnarray}
	\frac{\partial(\phi\rho_{\alpha}S_{\alpha})}{\partial t}=-\nabla\cdot(\rho_{\alpha}\textbf{u}_{\alpha})+q_{\alpha}~,\label{eq:mass}\\
	\textbf{u}_{\alpha}=-\frac{k_{r\alpha}}{\mu_{\alpha}}\textbf{K}(\nabla p_{\alpha}+\rho_{\alpha}g\nabla z)~,\label{eq:darcy}
\end{eqnarray}
subject to
\begin{eqnarray}
	S_w+S_o&=&1\\
	p_o-p_w&=&p_c
\end{eqnarray}
where $\alpha=w, o$ denotes the wetting and non-wetting phases, $S_{\alpha}=S_{\alpha}(\textbf{x}, t)$ is saturation, $p_{\alpha}=p_{\alpha}(\textbf{x}, t)$ is pressure.  $p_c=p_c(S_w)$ is capillary pressure, $\phi=\phi(p_o)$ is porosity, $\rho_{\alpha}=\rho_{\alpha}(p_{\alpha})$ is density, $\textbf{x}$ is spatial coordinate, $t$ is time, $\textbf{u}_{\alpha}$ is Darcy's velocity, $\textbf{K}=\textbf{K}(\textbf{x})$ is absolute permeability tensor, $k_{r\alpha}=k_{r\alpha}(S_{\alpha})$ is relative permeability, $\mu_{\alpha}=\mu_{\alpha}(p_{\alpha})$ is viscosity, $g$ is gravitational acceleration magnitude and $z$ is the vertical coordinate pointing upwards. All variables and meanings are summarized in Table \ref{tab:variables}. From Eqs.~\eqref{eq:mass} and \eqref{eq:darcy}, we have
\begin{equation}
	\frac{\partial(\phi\rho_{\alpha}S_{\alpha})}{\partial t}-\nabla\cdot(\frac{\rho_{\alpha}k_{r\alpha}}{\mu_{\alpha}}\textbf{K}(\nabla p_{\alpha}+\rho_{\alpha}g\nabla z))=q_{\alpha}~. \label{eq:2phase}
\end{equation}

In many practical scenarios of subsurface Darcy flows, the physical parameters of governing equations, such as permeability, porosity and relative permeability, are highly uncertain. The uncertain parameters are usually modelled as random variables in geostatistics \citep{caers2011modeling}. In this paper, permeability is assumed to be uncertain such that $K=K(\mathbf{x};\omega)$, where $\omega\in\Omega$ and $\Omega$ is the sample space. In consequence, we have $p_\alpha=p_\alpha(\mathbf{x},t; \omega)$ and $S_\alpha=S_\alpha(\mathbf{x},t; \omega)$ such that the governing equations \eqref{eq:2phase} become parametric or random PDEs.

\begin{table}[ht]
\centering
\caption{Variables and meanings.}
\label{tab:variables}
\begin{tabularx}{0.8\textwidth}{l X}
\toprule
Symbol & Description \\
\midrule
$S_\alpha$ & Saturation of phase $\alpha$ \\
$p_\alpha$ & Pressure of phase $\alpha$ \\
$p_c(S_w)$ & Capillary pressure \\
$\phi(p_o)$ & Porosity \\
$\mathbf{K}$ & Absolute permeability tensor \\
$k_{r\alpha}(S_\alpha)$ & Relative permeability of phase $\alpha$ \\
$\mu_\alpha(p_\alpha)$ & Dynamic viscosity of phase $\alpha$ \\
$\rho_\alpha(p_\alpha)$ & Density of phase $\alpha$ \\
$c_\alpha$ & Compressibility of phase $\alpha$ \\
$c_r$ & Rock compressibility \\
$q_\alpha$ & Source/sink term for phase $\alpha$ \\
\bottomrule
\end{tabularx}
\end{table}

\subsection{Operator Learning}
For the parametric governing equations \eqref{eq:2phase}, the objective of operator learning is to learn the mapping from $K(\mathbf{x};\omega)$ to $p_o(\mathbf{x},t; \omega)$ and $S_w(\mathbf{x},t; \omega)$ which are the primary solution variables. The governing equation can be compactly written as
\begin{equation}
	\mathcal{L}_{K(\omega)}u(\omega)=f
\end{equation}
where $\mathcal{L}_{K(\omega)}$ is a differential operator dependent on the random parameter field $K(\omega)$. $u(\omega)$ is the solution field denoting non-wetting-phase pressure and wetting-phase saturation such that $u(\omega)=[p_o(\omega), S_w(\omega)]$. This compact form represents a large class of governing equations with random parameter inputs arising in geoscience. Suppose the governing equation is consistent such that there is a unique $u$ for each realization of $K$, then a solution operator can be defined as
\begin{equation}
	\mathcal{G}: \mathcal{K}\rightarrow\mathcal{U},
\end{equation}
where $\mathcal{K}$ and $\mathcal{U}$ are suitable Banach spaces for $K\in\mathcal{K}$ and $u\in\mathcal{U}$. 

For problems involving random PDEs such as parameter inversion, uncertainty quantification and optimization, a large number of simulation runs are required to solve the governing PDE corresponding to each realization of $K$, which can be computationally very expensive. In operator learning, a model $\mathcal{G}_\theta$ with parameter $\theta$ is used to approximate $\mathcal{G}$ such that for a realization of $K$, the corresponding solution $u$ can be obtained efficiently as
\begin{equation}
	u\approx \mathcal{G}_\theta(K)~.
	\label{eq:g}
\end{equation} 

\section{Methodology}
\subsection{Numerical Simulation for Data Generation}
Operator learning for governing equations \eqref{eq:2phase} is based on data pairs $\{(K_i, u_i)\},i=1,...,N$ which can be generated by numerical simulation.
Adopting finite volume spatial discretization and fully implicit time integration, we have
\begin{equation}
	\frac{V_i}{\Delta t}(\phi_i^{n+1} S_{\alpha, i}^{n+1}\rho_{\alpha, i}^{n+1}-\phi_i^n S_{\alpha, i}^n\rho_{\alpha,i}^n)-\sum_j\rho_{\alpha, ij}^{n+1}\lambda_{\alpha, ij}^{n+1}T_{ij}(P_{\alpha, j}^{n+1}-P_{\alpha, i}^{n+1})-q_{\alpha, i}^{n+1}V_i=0~,\label{fv}
\end{equation}
for the $i$'th grid cell and $(n+1)$'th time step, where $ij$ denotes the face between $i$ and $j$'th cells, $P_\alpha=p_\alpha+\rho_\alpha gz$ is phase potential, $\rho_{\alpha, ij}^{n+1}$ is evaluated by arithmetic average, the mobility $\lambda_{\alpha, ij}^{n+1}=(k_{r\alpha}/\mu_{\alpha})_{ij}^{n+1}$ is evaluated by first-order upwind scheme, the transmissibility $T_{ij}=(\textbf{K}A/d)_{ij}$ is evaluated by harmonic average, $A$ is the area of shared face between neighbouring cells, and $d$ is the distance from cell centre to face centre. 
These are standard treatments for simulating compressible two-phase Darcy flows. 
The source terms are fixed or calculated using the Peaceman well model for grid cells containing injecting or producing wells. The functions of parameters can be found in Table~\ref{para1}.

State variables $S_w$ and $P_o$ are adopted as the primary variables. Hereafter, $S_w$ and $P_o$ are referred to as saturation and pressure for simplicity. As $\rho_\alpha$, $\phi$, $p_c$, $\mu_\alpha$ and $k_{r\alpha}$ are all functions of state variables, Eq.~\eqref{fv} leads to a nonlinear equation system

\begin{equation}
	\textbf{F}(\bm{\theta})=0 \label{nonlinear}
\end{equation} 
for $\bm{\theta}=[\textbf{S}_{w}^{n+1}, \textbf{P}_o^{n+1}]$ which is solved  by the Newton-Raphson method. For $\bm{\theta}_k$ of the $k$'th Newton-Raphson iteration,  $\Delta\bm{\theta}_k$ is obtained by solving a sparse linear equation system
\begin{equation}
	\textbf{J}(\bm{\theta}_k)\Delta\bm{\theta}_k=-\textbf{F}(\bm{\theta}_k)~, \label{n1}
\end{equation}
and $\bm{\theta}_k$ can be updated as 
\begin{equation}
	\bm{\theta}_{k+1}=\bm{\theta}_k+\Delta\bm{\theta}_k \label{n2}
\end{equation}
until convergence. 
It is worth mentioning that the assembling of $\textbf{J}(\bm{\theta}_k)$ and computation of the linear equation system $\textbf{J}(\bm{\theta}_k)\Delta\bm{\theta}_k=-\textbf{F}(\bm{\theta}_k)$ for $\Delta\bm{\theta}_k$ could be expensive using Python. We have employed Numba's just-in-time compilation with automatic parallelization and connected arrays to accelerate the assembly of the sparse Jacobian matrix. The LGMRES solver available in SciPy is used to solve the resulting linear equation system.

\begin{table}[!htp]
	\centering
	\begin{tabular}{|l|l|l|}
		\hline
		Parameter & Function & Hyperparameter (SI unit)\\
		\hline
		$\phi(p_o)$ & $\phi_r(1.0+c_r(p_o-p_r))$ & $\phi_r=0.1$, $p_r=2\times 10^7$, $c_r=10^{-11}$\\
		\hline
		$\rho_w(p_w)$ & $\rho_r10^{c_w(p_w-p_r)}$ & $\rho_r=10^3$, $p_r=2\times 10^7$, $c_w=10^{-9}$\\
		\hline
		$\rho_o(p_o)$ & $\rho_r10^{c_o(p_o-p_r)}$ & $\rho_r=850$, $p_r=2\times 10^7$, $c_o=5\times 10^{-8}$\\
		\hline
		$\mu_w(p_w)$ & $\mu_r{1+c_{\mu}(p_w-p_r)}$ & $\mu_r=10^{-3}$, $p_r=2\times 10^7$, $c_{\mu}=10^{-18}$\\
		\hline
		$\mu_o(p_o)$ & $\mu_r(1+c_{\mu}(p_o-p_r))$ & $\mu_r=2\times 10^{-3}$, $p_r=2\times 10^7$, $c_{\mu}=10^{-18}$\\
		\hline
		$k_{rw}(S_w)$ & $k_{rw0}\left(\frac{S_w-S_{wr}}{1.0-S_{wr}-S_{or}}\right)^{\beta}$ & $k_{rw0}=0.5$, $S_{wr}=0.2$, $S_{or}=0.2$, $\beta=2$\\
		\hline
		$k_{ro}(S_o)$ & $k_{ro0}\left(\frac{S_o-S_{or}}{1.0-S_{wr}-S_{or}}\right)^{\beta}$ & $k_{ro0}=1.0$, $S_{wr}=0.2$, $S_{or}=0.2$, $\beta=2$\\
		\hline
		$p_c(S_w)$ & $p_e$$\frac{1-S_w-S_{or}}{1-S_{or}-S_{wr}}$ & $p_e=10^4$, $S_{wr}=0.2$, $S_{or}=0.2$\\
		\hline
	\end{tabular}
	\caption{Parameter dependence on primary state variables. The hyperparameters are all in SI units.}
	\label{para1}
\end{table}

\subsection{Fourier Neural Operator for Cross-Resolution Surrogate Modelling}
The data pairs $\{(K_i, u_i)\},i=1,...,N_{samples}$ obtained by numerical simulation are finite-dimensional grid-based. The Fourier neural operator (FNO) is adopted as $\mathcal{G}_\theta$ to approximate the true inifinte-dimensional operator $\mathcal{G}$ based on finite-dimensional data such that cross-resolution surrogates can be obtained.
The input for FNO is the discrete permeability fields of a certain resolution, and the output is the pressure or wetting-phase saturation fields at the equivalent or higher resolution. Two FNO models, $G_{\theta}^{p}$ and $G_{\theta}^{S_w}$, are trained for pressure and saturation, respectively.

\begin{figure}[ht]
	\centering
	\vspace{0.0cm}
	\includegraphics[width=0.99\textwidth]{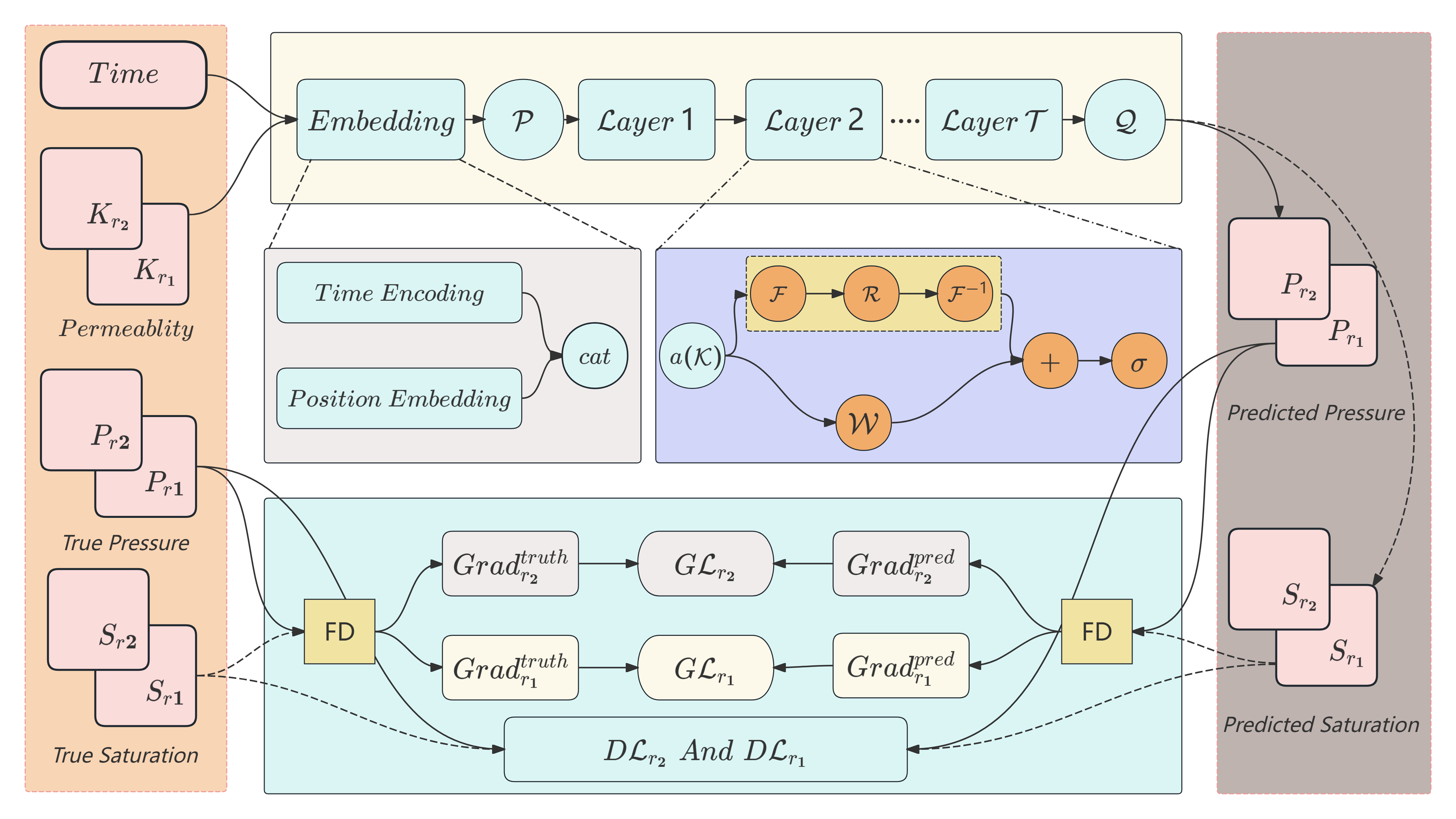}
	\caption{Neural Operator architecture for multiresolution Sobolev trained FNO for pressure or saturation. Permeability fields of low and medium resolutions, $K_{r_2}$ and $K_{r_1}$, and time step are used as the input. The input is augmented by the Embedding layer which embeds both temporal and spatial coordinates. A lifting operator $P$ maps the embeded input to a higher-dimensional feature space. The lifted features pass through multiple Fourier layers with nonlinear activation functions $\sigma$. Each Fourier layer consists of spectral convolution in the Fourier domain and a local linear transform. The projection operator $Q$ maps the learned features to the predicted pressure or saturation fields of resolutions $r_1$ and $r_2$ for training. The data losses $D\mathcal{L}_{r_2}$ and $D\mathcal{L}_{r_1}$ are computed using the MSE between predicted fields and ground truth. The gradient losses $G\mathcal{L}_{r_2}$ and $G\mathcal{L}_{r_1}$ are computed as the MSE between the gradient fields of the predicted pressure or saturation fields and those of the ground truth.}
	\label{fig:NN-structure}
\end{figure}

For the training data of a certain resolution level $n_z\times n_y\times n_x$, the input permeability tensor is of dimension $n_b\times 1\times n_z\times n_y\times n_x$, and the output is of dimension $n_b\times 1\times n_z\times n_y\times n_x$, where $n_b$ is the batch size. The input is one permeability field and time step; the output is one pressure field for $G_{\theta}^{p}$ and saturation field for $G_{\theta}^{S_w}$. Both inputs and ouputs are normalized to be within [0,1]. 

The FNO-based $G_\theta$ is written as
\begin{equation}
G_\theta = Q \circ (W_L + K_L) \circ \cdots \circ \sigma(W_1 + K_1) \circ P~.
\end{equation}
The FNO structure (see Fig.~\ref{fig:NN-structure}) stacks four Fourier layers, each with a uniform hidden width of 128. The network incorporates positional embedding to explicitly map three-dimensional spatial coordinates into a high-dimensional feature space, providing the model with continuous spatial awareness and facilitating the learning of coordinate-dependent solution patterns. Within each Fourier layer, a linear transformation is applied to the input in the frequency domain via the Fast Fourier Transform, and the frequency components are truncated. The truncation modes along the z, y, and x axes are set to 5, 8, and 8, respectively.

AdamW with learning rate $0.005$ and a cosine-annealing schedule with gradient clipping are adopted for optimization. 
FNO is trained using labeled samples of different lower resolutions and generalized to samples of higher resolution.

\subsection{Multiresolution Sobolev Learning}
In the current study, a multiresolution Sobolev training method is proposed to enhance the cross-resolution prediction accuracy using spatial gradients as regularizers in the loss function. 
For FNO outputs, we only have solution values at fixed grid cells. Therefore, gradients are computed by the finite difference scheme where the discretization intervals are consistent with the grid of numerical simulation. This implies that the discretization intervals for approximating gradients are resolution-dependent.
Let the grid cell is labeled as $i,j,k$ on a Cartetian grid of dimension $n_x \times n_y \times n_z$, the first-order gradients in the $x$ direction are computed as

\begin{equation}
\partial_x {u}_{i,j,k}\approx
\begin{cases}
\dfrac{{u}_{i+1,j,k}-u_{i-1,j,k}}{2\Delta x}, & 2\le x\le n_x-1,\\[6pt]
\dfrac{u_{2,j,k}-u_{1,j,k}}{\Delta x}, & i=1,\\[6pt]
\dfrac{u_{n_x,j,k}-u_{n_x-1,j,k}}{\Delta x}, & i=n_x,
\end{cases}
\end{equation}

First-order gradients in $y$ and $z$ directions can be obtained similarly. Further, second-order gradients can be computed analogously based on first-order gradients. 
The total loss function based on the mean squared error (MSE) between prediction and labeled data for $u\in \{P, S_w\}$ is
\begin{equation}
	\mathcal{L}^{\text{total}}_{u} = \sum_{p=1}^{N_p}w^p\mathcal{L}^{p}_{u}~,
\end{equation}
where $\mathcal{L}^{p}_{u}$ denotes loss term related to the $p$'th-order gradients with weighting coefficient $w^p$. The order of $\mathcal{L}^{\text{total}}_{u}$ refers to the highest order of gradients in the loss, i.e., $N_p$. We compute until the second-order gradients such that $N_p=2$. The zero-order term with $p=0$ is the original data-driven loss functional that can be defined as the expectation of squared $L_2$ distance between predicted $u$ and the ground truth for single-resolution cases as
\begin{equation}
	\hat{\mathcal{L}}_{u}^{0}=\mathbb{E}_\omega[|u-u^*|_2^2]
\end{equation}
where the dimension of $u$ is $N_u=N_z\times N_y\times N_x$. However, $|u-u^*|_2^2$ is affected by $n_u$ representing the resolution of $u$. For cross-resolution operator learning, the zero-order loss functional is defined as the expectation of mean squared error (MSE) between predicted $u$ and the ground truth as
\begin{equation}
	\mathcal{L}_{u}^{0}=\mathbb{E}_\omega\left[\frac{1}{N_u}|u-u^*|_2^2\right]~,
\end{equation}
which can be approximated as
\begin{equation}
	\mathcal{L}_{u}^{0}\approx \frac{1}{N_{samples}}\sum_{r=1}^{N_{samples}}\left[\frac{1}{N_u}|u_r-u_r^*|_2^2\right]
    \label{zeroorder}
\end{equation}
on $N_{samples}$ training samples denoted by $u_r$. Minimizing $\mathcal{L}_{u}^{0}$ is equivalent to minimizing $\hat{\mathcal{L}}_{u}^{0}$ for fixed $N_u$, but adopting $\mathcal{L}_{u}^{0}$ can eliminate the effect of resolutions on the loss function for training samples with mixed resolutions, which is key for learning cross-resolution surrogates.
Similarly, the first-order gradient loss $\mathcal{L}_{u}^{1}$ is approximated as
\begin{equation}
	\mathcal{L}_{u}^{1}\approx \frac{1}{N_{samples}}\sum_{r=1}^{N_{samples}}\left[\frac{1}{N_u}(|\partial_xu_r-\partial_xu_r^*|_2^2+|\partial_yu_r-\partial_yu_r^*|_2^2+|\partial_zu_r-\partial_zu_r^*|_2^2)\right]~.
\end{equation}
Further, the second-order gradient loss $\mathcal{L}_{u}^{2}$ is approximated as
\begin{equation}
	\begin{split}
	\mathcal{L}_{u}^{2}\approx \frac{1}{N_{samples}}\sum_{r=1}^{N_{samples}}\left[\frac{1}{N_u}\left(|\partial_{xx}u_r-\partial_{xx}u_r^*|_2^2+|\partial_{xy}u_r-\partial_{xy}u_r^*|_2^2+|\partial_{xz}u_r-\partial_{xz}u_r^*|_2^2\right.\right.\\
	|\partial_{yx}u_r-\partial_{yx}u_r^*|_2^2+|\partial_{yy}u_r-\partial_{yy}u_r^*|_2^2+|\partial_{yz}u_r-\partial_{yz}u_r^*|_2^2\\
	\left.\left.|\partial_{zx}u_r-\partial_{zx}u_r^*|_2^2+|\partial_{zy}u_r-\partial_{zy}u_r^*|_2^2+|\partial_{zz}u_r-\partial_{zz}u_r^*|_2^2\right)\right]~.
	\end{split}
\end{equation}

\subsection{Convergence of Sobolev Training}
The effectiveness of Sobolev training for acclelerating the convergence of NN under certain assumptions have been proved in \citep{cho2025sobolev}. The main result is presented here concisely for providing a reference of theoretical guarantee although the assumptions are not valid for FNO. Generalizing the result to the complexity of FNO is highly non-trivial. The current study establishes a multiresolution Sobolev training algorithm for FNO and the effectiveness for enhancing convergence as well as cross-resolution generalization are valided by numerical experiments. The theoretical convergence analysis in \citet{cho2025sobolev} is based on the following assumptions.

\begin{assumption} \label{ass:main}
Let $\{(v_k(x),u_k(x))\}_{k=1}^N$ be a pair of input and output real-valued functions, $\omega$ denotes NN parameters. The following conditions hold.
\begin{enumerate}
\item The query points $\{x_j\}_{j=1}^J \subset \mathbb{R}^n$ are i.i.d. sampled from a spherical Gaussian distribution $\mathcal{N}(0, I_n)$.
\item The output function $u_k(x)$ is approximated by an integral operator
\[
u_k(x) \approx \frac{1}{J}\sum_{\ell=1}^J \kappa(x,y_\ell) v_k(y_\ell),
\]
with a kernel $\kappa(x,y)$ that admits the separable form
\[
\kappa(x,y)=\sigma(w^\top x)\,\sigma(w^\top y),
\]
where $\sigma(\cdot)=\mathrm{ReLU}(\cdot)$.
\item The input functions are
\[
v_k(x)=\mu_k\,\sigma((w^*)^\top x),\qquad \mu_k\in\mathbb{R},\; w^*\in\mathbb{R}^n,
\]
and the output functions predicted by the network are
\[
u_k(x;w)=\frac{1}{J}\sum_{\ell=1}^J \mu_k\,\sigma(w^\top x)\,\sigma(w^\top y_\ell)\,\sigma((w^*)^\top y_\ell).
\]
The ground‑truth outputs are $u_k(x;w^*)$.
\item The parameters evolve according to the gradient flows
\[
\frac{dw_{L^2}}{dt}=-\mathbb{E}_X\Bigl[\frac{\partial\mathcal{L}_{L^2}}{\partial w}\Bigr],\qquad
\frac{dw_{\mathrm{Sob}}}{dt}=-\mathbb{E}_X\Bigl[\frac{\partial\mathcal{L}_{\mathrm{Sob}}}{\partial w}\Bigr],
\]
with
\[
\mathcal{L}_{L^2}=\frac1N\sum_{k=1}^N\bigl\|u_k(\cdot;w)-u_k(\cdot;w^*)\bigr\|_{L^2}^2,\qquad
\mathcal{L}_{\mathrm{Sob}}=\mathcal{L}_{L^2}+\mathcal{L}_{\mathrm{Der}},
\]
\[
\mathcal{L}_{\mathrm{Der}}=\frac1N\sum_{k=1}^N\bigl\|D_x u_k(\cdot;w)-D_x u_k(\cdot;w^*)\bigr\|_{L^2}^2.
\]
\item The initial weight vector satisfies $\|w_0-w^*\|<\|w^*\|$.
\end{enumerate}
\end{assumption}

The theoretical guarantee for Sobolev training is by the following theorem in \citet{cho2025sobolev} where the complete proof can be found.
\begin{theorem}[Convergence of Sobolev training] \label{thm:main}
Under Assumption~\ref{ass:main}, we have
\[
\frac{d}{dt}\|w_{L^2}-w^*\|^2 < 0
\]
whenever $\|w_0-w^*\|<\|w^*\|$, and $w_{L^2}(t)\to w^*$ as $t\to\infty$.
For the Sobolev loss we have
\[
\frac{d}{dt}\|w_{\mathrm{Sob}}-w^*\|^2 \;\le\; \frac{d}{dt}\|w_{L^2}-w^*\|^2 \;<\; 0.
\]
Moreover, if $w$ and $w^*$ are not parallel (i.e. the angle $\theta(w,w^*)\in(0,\pi)$), then the inequality is strict
\[
\frac{d}{dt}\|w_{\mathrm{Sob}}-w^*\|^2 \;<\; \frac{d}{dt}\|w_{L^2}-w^*\|^2 \;<\; 0.
\]
Hence incorporating derivative information accelerates the convergence rate.
\end{theorem}

\section{Numerical Experiments}
\label{sec:num_exp}
In this section, numerical experiments are conducted to validate the accuracy of FNO-based surrogate with multiresolution Sobolev learning for cross-resolution prediction of 3D pressure and saturation fields in two-phase compressible Darcy flow examples.

\subsection{Multiresolution Data Generation}
For data generation, 3D permeability fields of resolutions $11\times11\times6$, $21 \times 21\times6$ and $41\times41\times6$ are generated. The horizontal logrithm-permeability fields are generated by sequential Gaussian simulation using the Stanford Geostatistical Modelling Software (SGeMS) \citep{remy2009applied}. The permeability tensor is assumed to be diagonal such that the horizontal permeabilities are identical while the vertial permeability is one tenth of the horizontal conterpart. The size of the reservoir model is 300 m$\times$300 m$\times$30m.

\begin{figure}[ht]
	\centering
	\vspace{0.0cm}
	\includegraphics[width=0.99\textwidth]{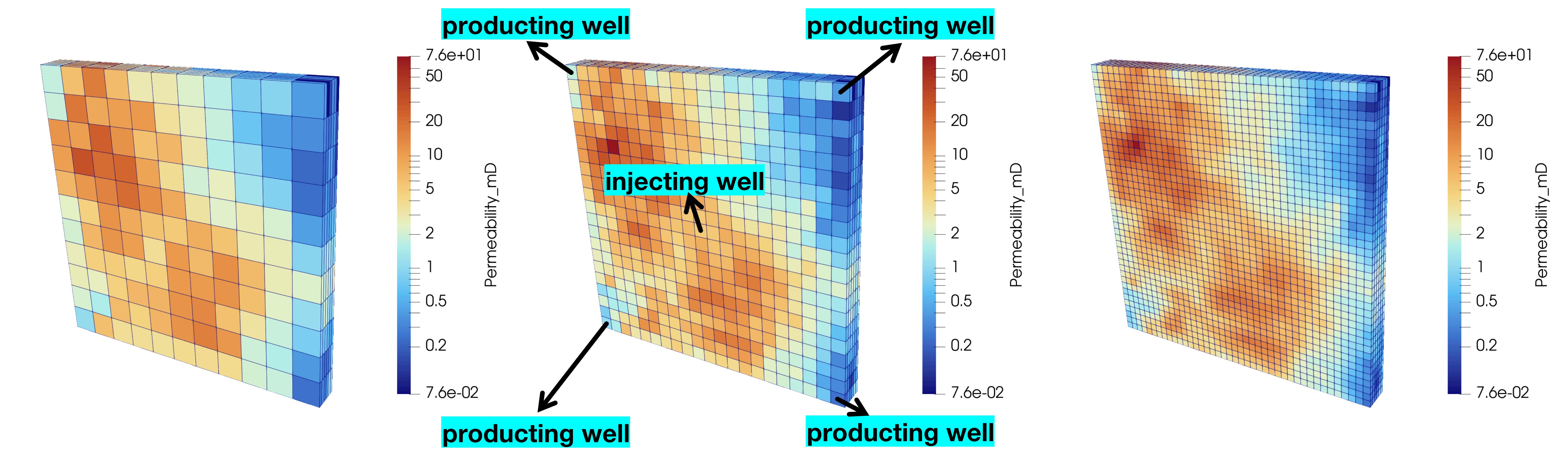}
	\caption{Permeability fields at low, medium and high resolutions.}
	\label{fig:perm_differ_resolution}
\end{figure}

\begin{figure}[ht]
	\centering
	\vspace{0.0cm}
	\includegraphics[width=0.99\textwidth]{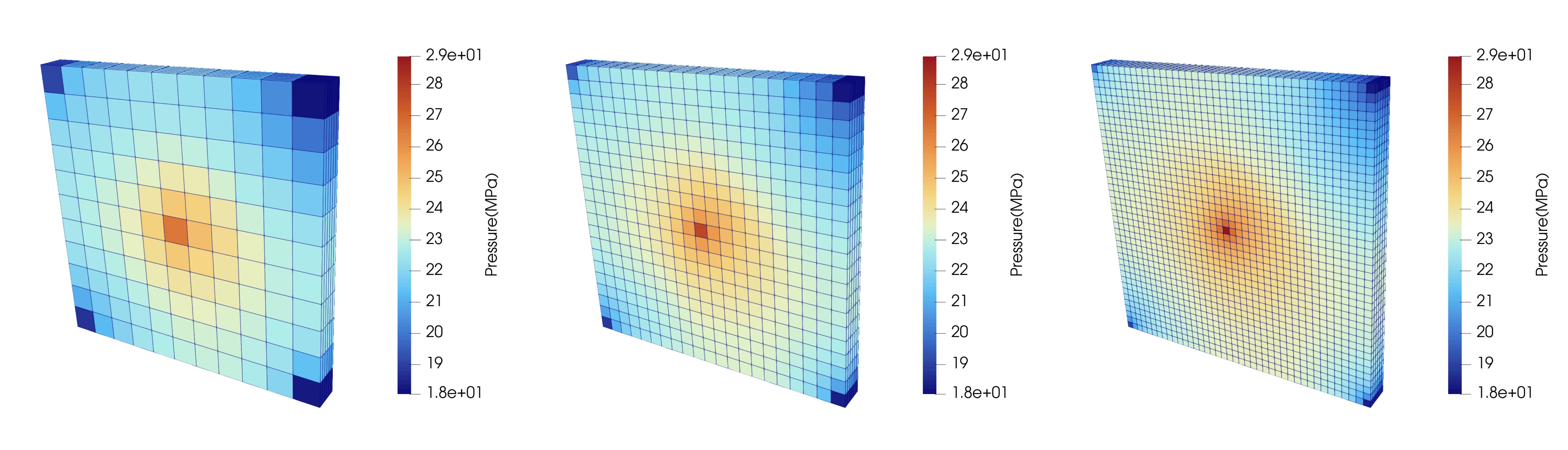}
	\caption{Pressure fields at low, medium and high resolutions.}
	\label{fig:press_differ_resolution}
\end{figure}

\begin{figure}[ht]
	\centering
	\vspace{0.0cm}
	\includegraphics[width=0.99\textwidth]{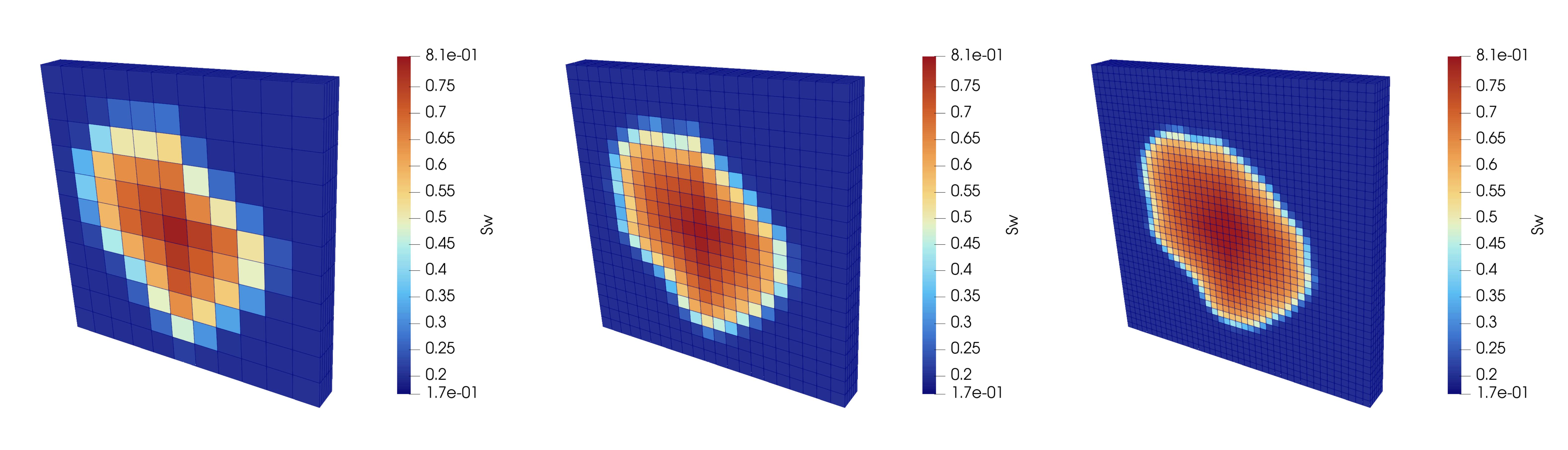}
	\caption{Saturation fields at low, medium and high resolutions.}
	\label{fig:sw_differ_resolution}
\end{figure}

The transient pressure and wetting-phase saturation fields are generated by numerical simulation as the labeled data. For simulation, there is an injecting well at the center of the reservoir, and four producting wells at the corners (see Fig.~\ref{fig:press_differ_resolution}). The injector is at fixed flow rate while the producers are at fixed bottom-hole pressure. The initial conditions of simulation are pressure at 20 MPa and wetting-phase saturation at 0.2 everywhere in the reservoir. The injecting well is at constant flow rate of 10 $m^3$ per day while the producing wells are at constant bottom-hole pressure of 18 MPa. Representative permeability fields at low, medium and high resolutions are shown in Fig.~\ref{fig:perm_differ_resolution}. The pressure and saturation fields at end of simulation on these samples are presented in Fig.~\ref{fig:press_differ_resolution} and Fig.~\ref{fig:sw_differ_resolution}.

\subsection{Training Processes}
\begin{figure}[ht]
	\centering
	\vspace{0.0cm}
	\includegraphics[width=0.8\textwidth]{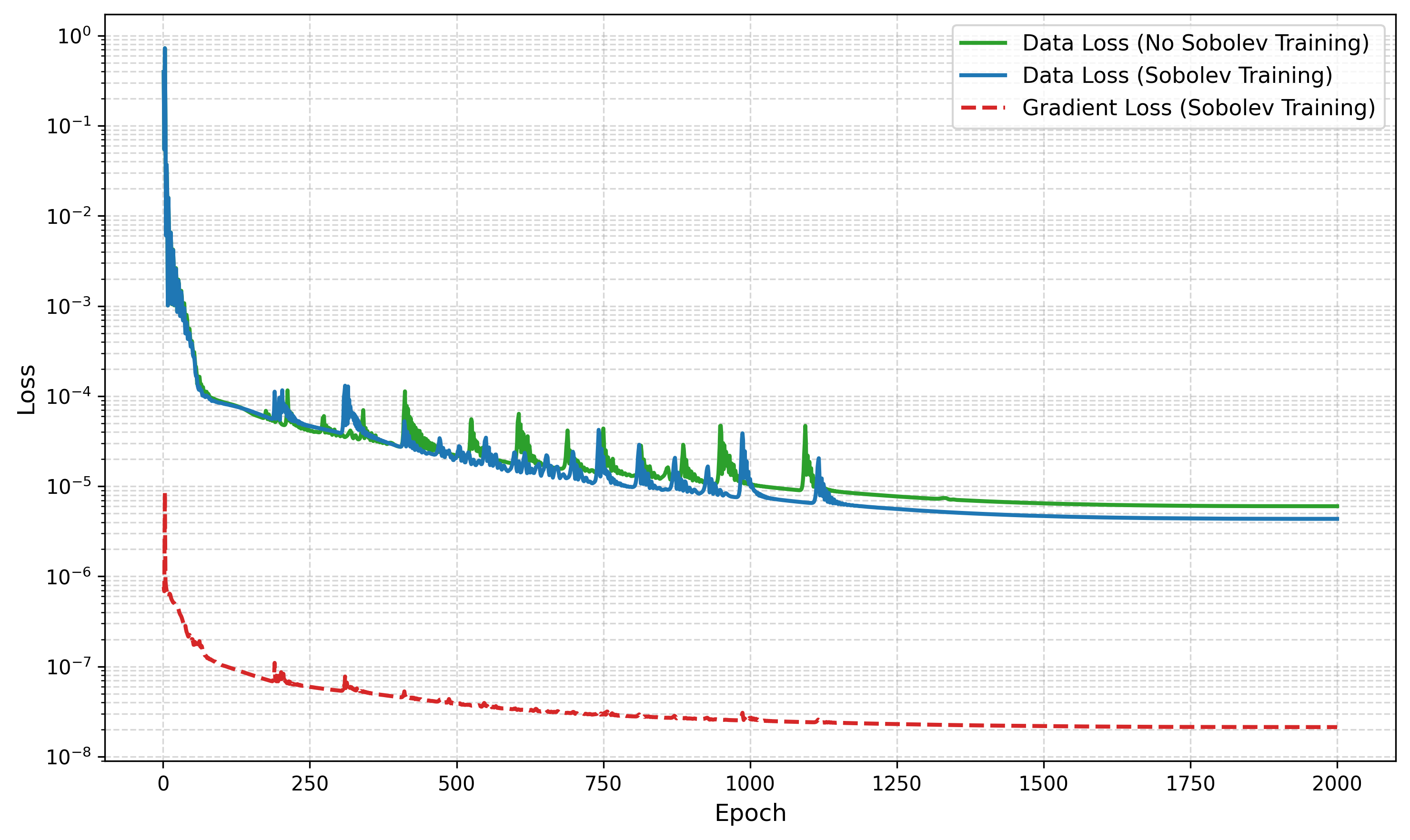}
	\caption{Comparison of training losses in conventional training without gradients and Sobolev training incorporating first-order gradients into the loss function for surrogate modelling of the pressure field based on FNO. }
	\label{fig:Press_Loss_Comparison}
\end{figure}

\begin{figure}[ht]
	\centering
	\vspace{0.0cm}
	\includegraphics[width=0.8\textwidth]{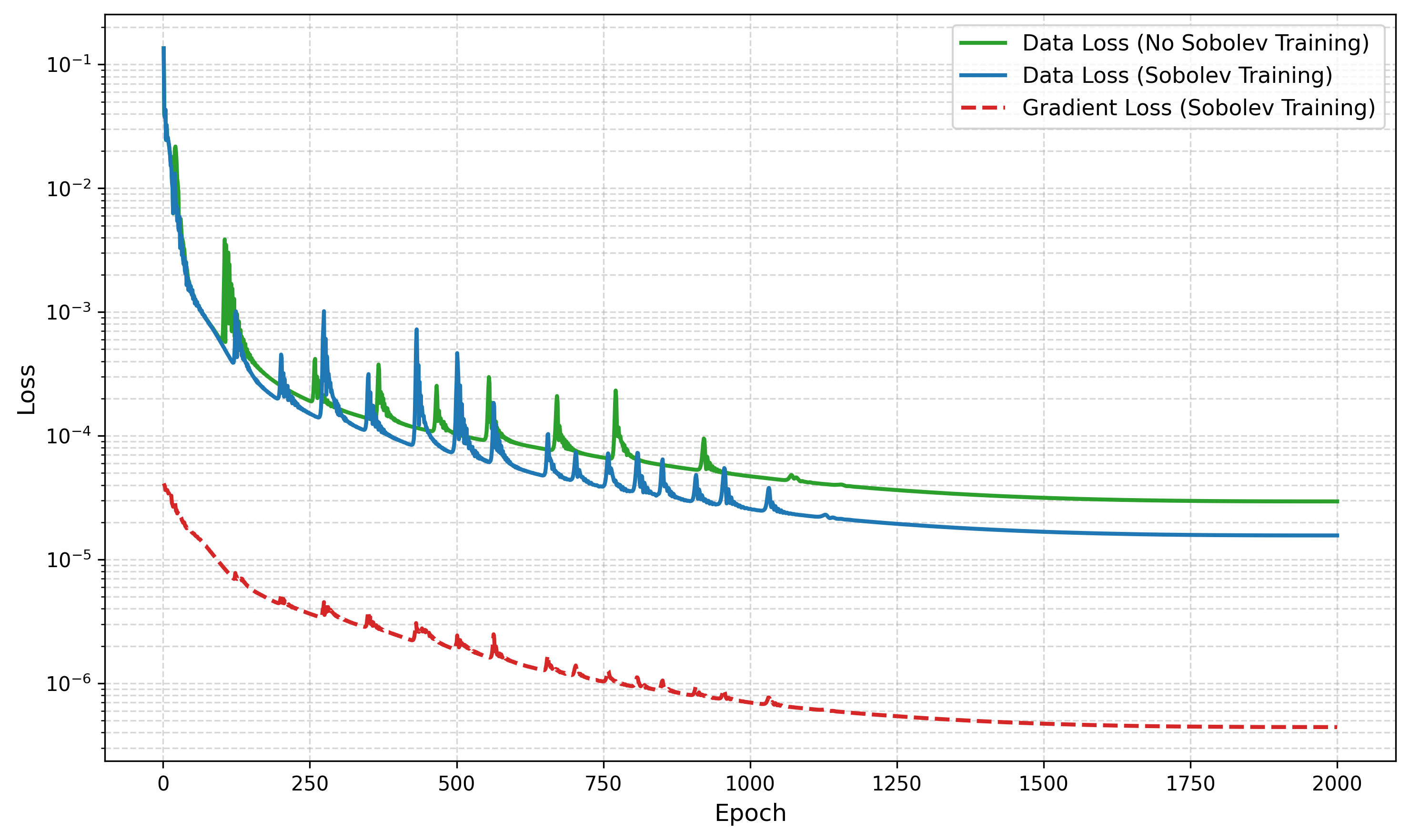}
	\caption{Comparison of training losses in conventional training without gradients and Sobolev training incorporating first-order gradients into the loss function for surrogate modelling of the saturation field based on FNO.}
	\label{fig:Sw_Loss_Comparison}
\end{figure}

The training process of FNO for pressure and saturation fields are compared in Fig.~\ref{fig:Press_Loss_Comparison} and Fig.~\ref{fig:Sw_Loss_Comparison}, where data loss refers to the zero-order loss term $\mathcal{L}_u^0$ while gradient loss is the first-order loss term $\mathcal{L}_u^1$. For Fig.~\ref{fig:Press_Loss_Comparison}, the training set contains 200 labelled samples of resolution $21\times21\times6$ and 100 labelled samples of resolution $11\times11\times6$ for the pressure field. For Fig.~\ref{fig:Sw_Loss_Comparison}, the training set contains 200 labelled samples of resolution $21\times21\times6$ and 100 labelled samples of resolution $11\times11\times6$ for saturation. The comparison shows that incorporating gradients into the loss function helps accelerate the convergence of $\mathcal{L}_u^0$ for training FNO with multiresolution labelled data.

\subsection{Effect of Cross-Resolution Training set}
Trained on lower resolutions, the FNO-based surrogate can be used to predict the pressure and saturation fields of higher resolution. 
We demonstrate that by using cross-resolution training data, the cross-resolution prediction accuracy for samples of higher resolution can be considerably increased. In the numerical experiment, the FNO-based surrogate is trained on labelled data of resolutions $21\times21\times6$ and $11\times11\times6$, and  is subsequently used to predict the pressure and saturation fields of resolution $41\times41\times6$ at end of simulation.

Table \ref{tab:FNO_press_cross1} summaries the prediction accuracy of surrogates using different training sets for the pressure field. Table \ref{tab:FNO_sw_cross1} summaries the corresponding results for the saturation field. The prediction accuracy is quantified by the coefficient of determination $R^2$ and mean absolute error (MAE) to guarantee that the surrogate models are both accurate and more effective than using the mean value of labelled data.
It is observed from the tables that by including some low-resolution samples that that cheaper to generate can yield much higher cross-resolution prediction accuracy.

\begin{table}[H]
\centering
\begin{tabular}{cccc}
\toprule
Resolution-1 & Resolution-2 & $R^2$ & MAE (MPa)\\
\midrule

150 samples &   0 samples & -4.080525 & 0.5127702 \\
150 samples &  50 samples &  0.321672 & 0.15894451 \\
200 samples &   0 samples & -3.930501 & 0.5049159 \\
200 samples & 100 samples &  0.411321 & 0.14184923 \\

\bottomrule                            
\end{tabular}
\caption{Prediction accuracy using different training sets for surrogate modelling based on FNO for the pressure field. The first two columns are number of samples of resolution-1 ($21\times21\times6$) and resolution-2 ($11\times11\times6$). The prediction accuracy is quantified on samples of resolution $41\times41\times6$. The MAE is averaged on both cells and samples.}
\label{tab:FNO_press_cross1}
\end{table}

\begin{table}[H]
\centering
\begin{tabular}{cccc}
\toprule
Resolution-1 & Resolution-2 & $R^2$ & MAE (percent)\\
\midrule

150 samples &   0 samples & 0.519655 & 1.605554 \\
150 samples &  50 samples & 0.747141 & 1.157219 \\
200 samples &   0 samples & 0.608308 & 1.512836 \\
200 samples & 100 samples & 0.757315 & 1.140762 \\

\bottomrule
\end{tabular}
\caption{Prediction accuracy using different training sets for surrogate modelling based on FNO for the saturation field. The first two columns are number of samples of resolution-1 ($21\times21\times6$) and resolution-2 ($11\times11\times6$). The prediction accuracy is quantified on samples of resolution $41\times41\times6$. The MAE is averaged on both cells and samples.}
\label{tab:FNO_sw_cross1}
\end{table}

\begin{figure}[!htb]
	\centering
	\vspace{0.0cm}
	\includegraphics[width=0.99\textwidth]{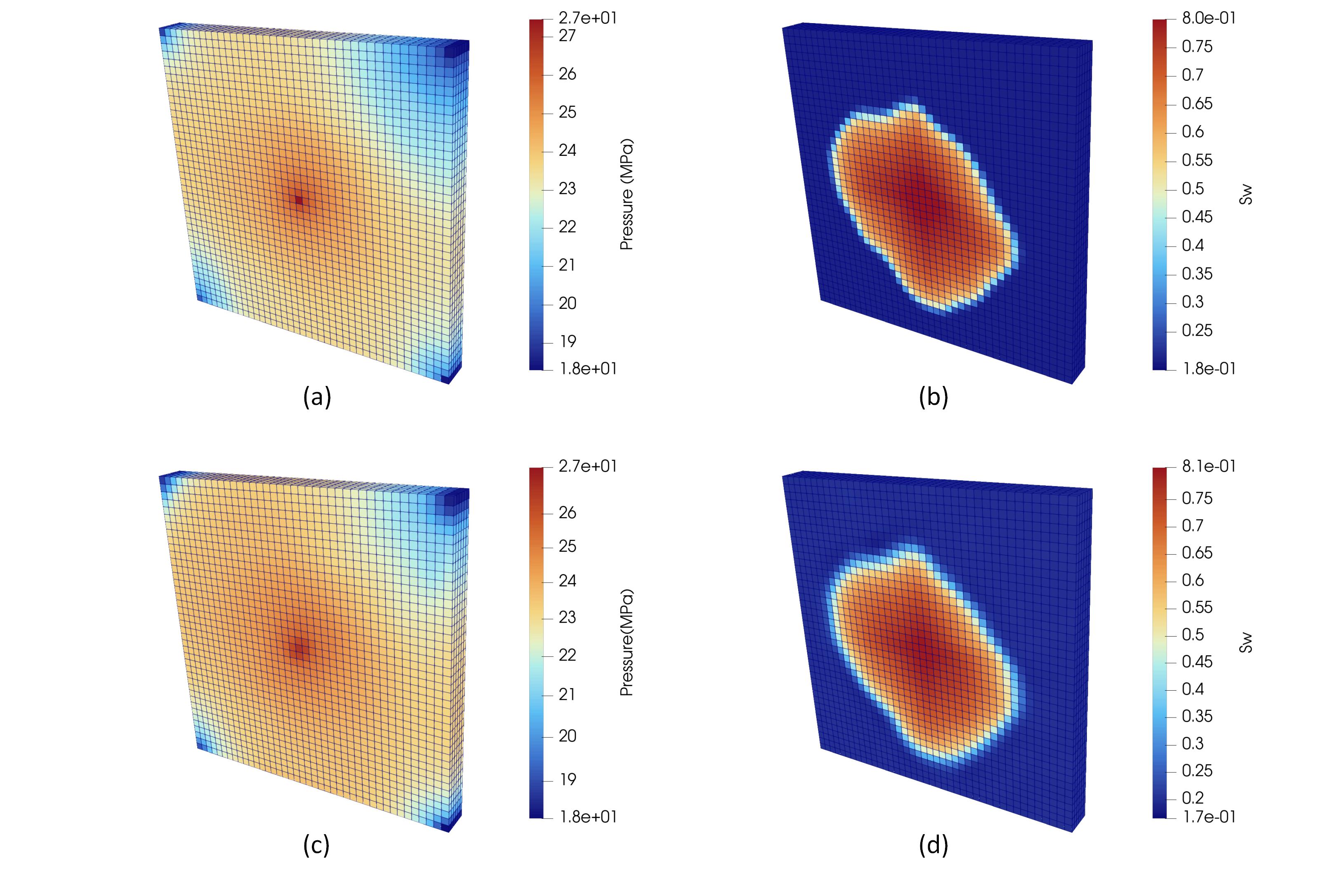}
	\caption{Comparison of true and predicted solution fields:
True pressure (a) and saturation (b) fields; Predicted pressure (c) and saturation fields (d). The training set consists of 150 $21\times21\times6$ samples and 50 $11\times11\times6$ samples. The presented pressure and saturation field of resolution $41\times41\times6$ corresponds to a sample of medium prediction accuracy in the test set.}
	\label{fig:figcomparetrue.jpg}
\end{figure}

\begin{figure}[!htb]
	\centering
	\vspace{0.0cm}
	\includegraphics[width=0.99\textwidth]{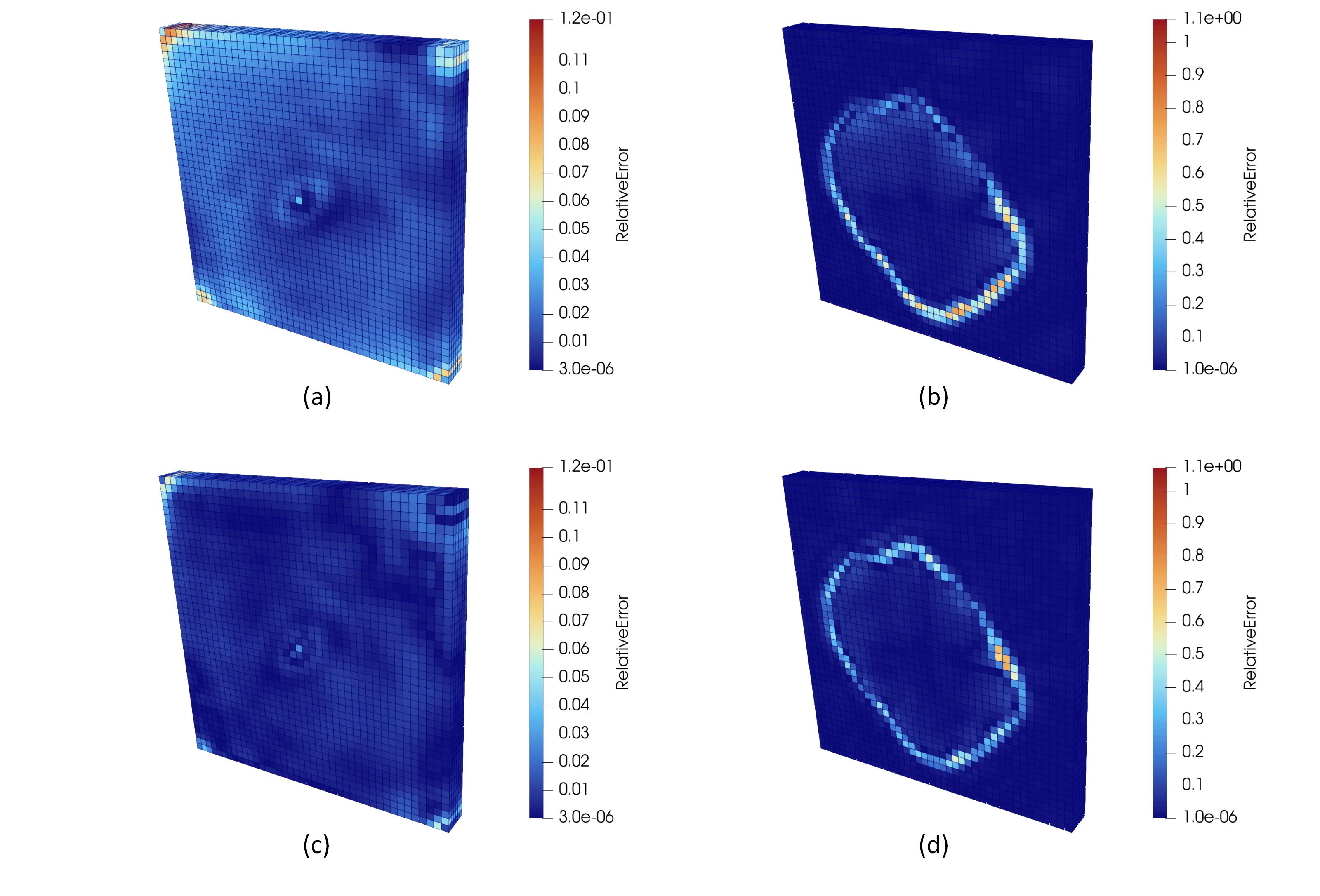}
	\caption{The heat map of relative error for predicted pressure and saturation fields on a sample of medium accuracy in the test set.
Pressure field using 200 $21\times21\times6$ samples (a); saturation field using 200 $21\times21\times6$ samples (b); pressure field using 150 $21\times21\times6$ samples and 50 $11\times11\times6$ samples (c); saturation field using 150 $21\times21\times6$ samples and 50 $11\times11\times6$ samples (d).}
	\label{fig:comparingerrorheatmap.jpg}
\end{figure}

\begin{figure}[!htb]
	\centering
	\vspace{0.0cm}
	\includegraphics[width=0.8\textwidth]{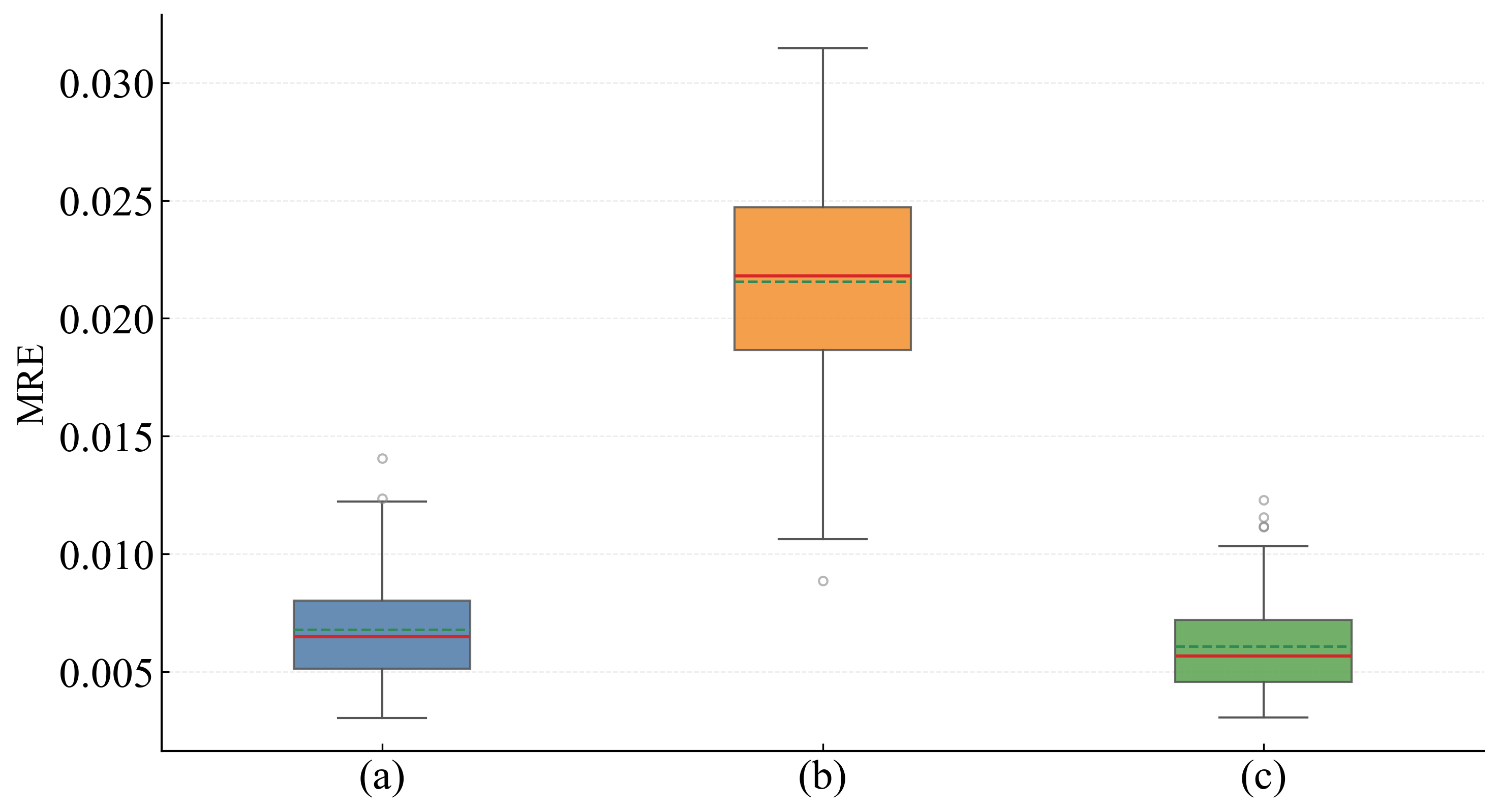}
	\caption{Box-plots of prediction accuracy for pressure using MRE comparing different training set compositions. The MRE is obtained by averaging relative errors on cells. The box plots show the uncertainty of MRE for different samples. The comparison is for training sets consisting of 150 $21\times21\times6$ samples and 50 $11\times11\times6$ samples (a); 200 $21\times21\times6$ samples (b); 200 $21\times21\times6$ samples and 100 $11\times11\times6$ samples (c).}
	\label{fig:sample_differ_pressure_of_box_plots_in_MRE_comparing}
\end{figure}

\begin{figure}[!htb]
	\centering
	\vspace{0.0cm}
	\includegraphics[width=0.8\textwidth]{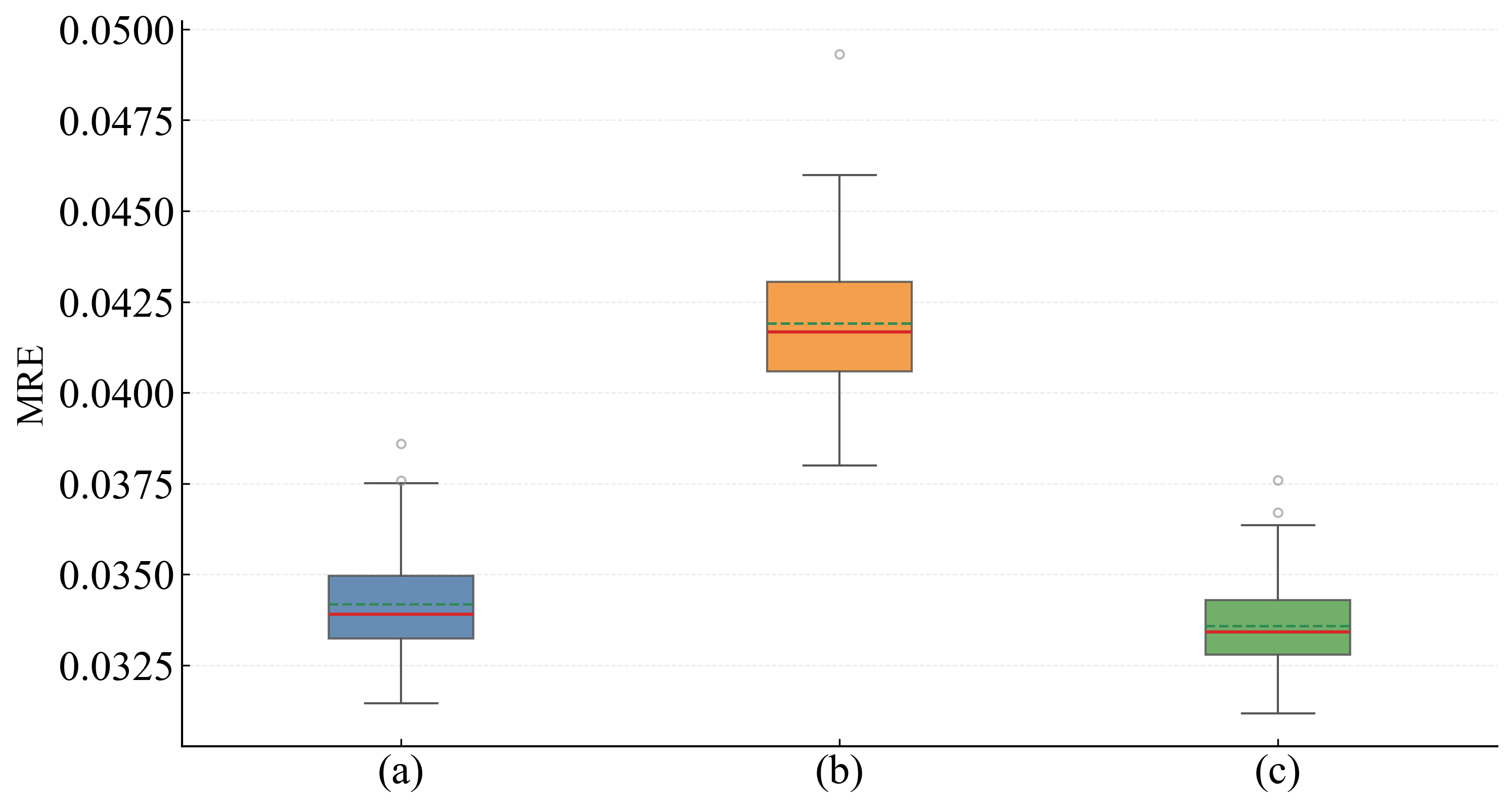}
	\caption{Box-plots of prediction accuracy for saturation using MRE comparing different training set compositions. The MRE is obtained by averaging relative errors on cells. The box plots show the uncertainty of MRE for different samples. The comparison is for training sets consisting of 150 $21\times21\times6$ samples and 50 $11\times11\times6$ samples (a); 200 $21\times21\times6$ samples (b); 200 $21\times21\times6$ samples and 100 $11\times11\times6$ samples (c).}
	\label{fig:sample_differ_saturation_of_box_plots_in_MRE_comparing}
\end{figure}

The true and predicted fields at end of simulation are presented in Fig.~\ref{fig:figcomparetrue.jpg}. The pressure fields are visually very close, while there is noticeable difference at the water fronts in the saturation fields. The MRE fields of FNO trained using single and multiresolution training sets, keeping the total number of training samples identical, are shown in Fig.~\ref{fig:comparingerrorheatmap.jpg}. While cheaper to be generated, the multiresolution training sets yield lower MRE, especially at the water front in the saturation field. 
To account for the uncertainty of prediction accuracy, box-plots of MRE averaged on grid cells for different samples in the test set are presented in Fig.~\ref{fig:sample_differ_pressure_of_box_plots_in_MRE_comparing} and Fig.~\ref{fig:sample_differ_saturation_of_box_plots_in_MRE_comparing}, where the effect of using multi-resolution training samples are notable.

Further, we investigate the effect of two different settings for the cross-resolution training sets. In setting 1, the lower resolution-2 $11\times11\times6$ samples are generated as new samples. In setting 2, the lower resolution-2 samples are generated by interpolating medium resolution-1 $21\times21\times6$ samples. In both settings, the higher-resolution $41\times41\times6$ samples in the test set are generated as new samples. The results for pressure and saturation are summarised in Tables \ref{tab:FNO_press_cross2} and \ref{tab:FNO_sw_cross2}. The results show that using newly generated lower resolution samples can yield higher prediction accuracy than using interpolated lower resolution samples.

\begin{table}[!htb]
\centering
\begin{tabular}{cccccc}
\toprule
Training set & Resolution-1 & Resolution-2 & Resolution-2 Generation&$R^2$ & MAE (MPa)\\
\midrule

1 & 100 samples &  50 samples & new & 0.289525 & 0.15984024 \\
2 & 150 samples &  50 samples & new & 0.321672 & 0.15894451 \\
3 & 200 samples & 100 samples & new & 0.411321 & 0.14184923 \\
4 & 100 samples &  50 samples & Interpolation & 0.191400 & 0.17595704 \\
5 & 150 samples &  50 samples & Interpolation & 0.331139 & 0.15438825 \\
6 & 200 samples & 100 samples & Interpolation & 0.373178 & 0.14918774 \\

\bottomrule
\end{tabular}
\caption{Prediction accuracy comparing different training set compositions for the pressure field. The first two columns are number of samples of resolution-1 ($21\times21\times6$) and resolution-2 ($11\times11\times6$). The prediction accuracy is quantified on data of resolution $41\times41\times6$. The MAE is averaged on both cells and samples.}
\label{tab:FNO_press_cross2}
\end{table}

\begin{table}[!htb]
\centering
\begin{tabular}{cccccc}
\toprule
Training set & Resolution-1 & Resolution-2 & Resolution-2 Generation&$R^2$ & MAE (percent)\\
\midrule

1 & 100 samples &  50 samples & new & 0.738507 & 1.154046 \\
2 & 150 samples &  50 samples & new & 0.747141 & 1.157219 \\
3 & 200 samples & 100 samples & new & 0.757315 & 1.140762 \\
4 & 100 samples &  50 samples & Interpolation & 0.734825 & 1.15636 \\
5 & 150 samples &  50 samples & Interpolation & 0.747899 & 1.152536 \\
6 & 200 samples & 100 samples & Interpolation & 0.753650 & 1.143951 \\

\bottomrule
\end{tabular}
\caption{Prediction accuracy comparing different training set compositions for the saturation field. The first two columns are number of samples of resolution-1 ($21\times21\times6$) and resolution-2 ($11\times11\times6$). The prediction accuracy is quantified on data of resolution $41\times41\times6$. The MAE is averaged on both cells and samples.}
\label{tab:FNO_sw_cross2}
\end{table}

\subsection{Effect of Sobolev Training}
\begin{figure}[!htb]
	\centering
	\vspace{0.0cm}
	\includegraphics[width=0.7\textwidth]{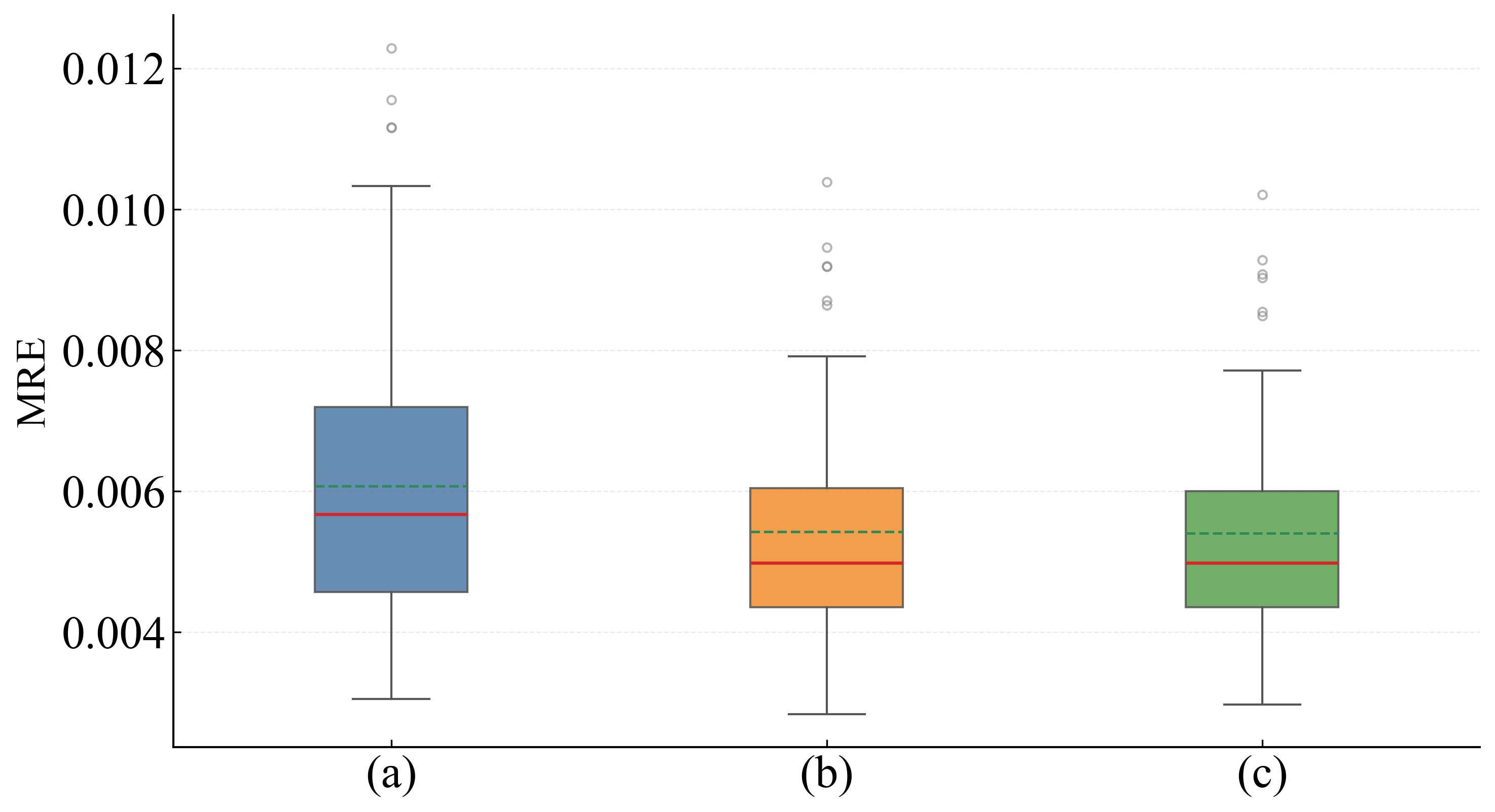}
	\caption{Prediction MRE for $41\times41\times6$ pressure fields using zero-order (a), first‑order (b) and second‑order (c) loss function. The training set consists of 200 $21\times21\times6$ and 100 $11\times11\times6$ samples.}
	\label{fig:press41_compare_200_100}
\end{figure}

\begin{figure}[!htb]
	\centering
	\vspace{0.0cm}
	\includegraphics[width=0.7\textwidth]{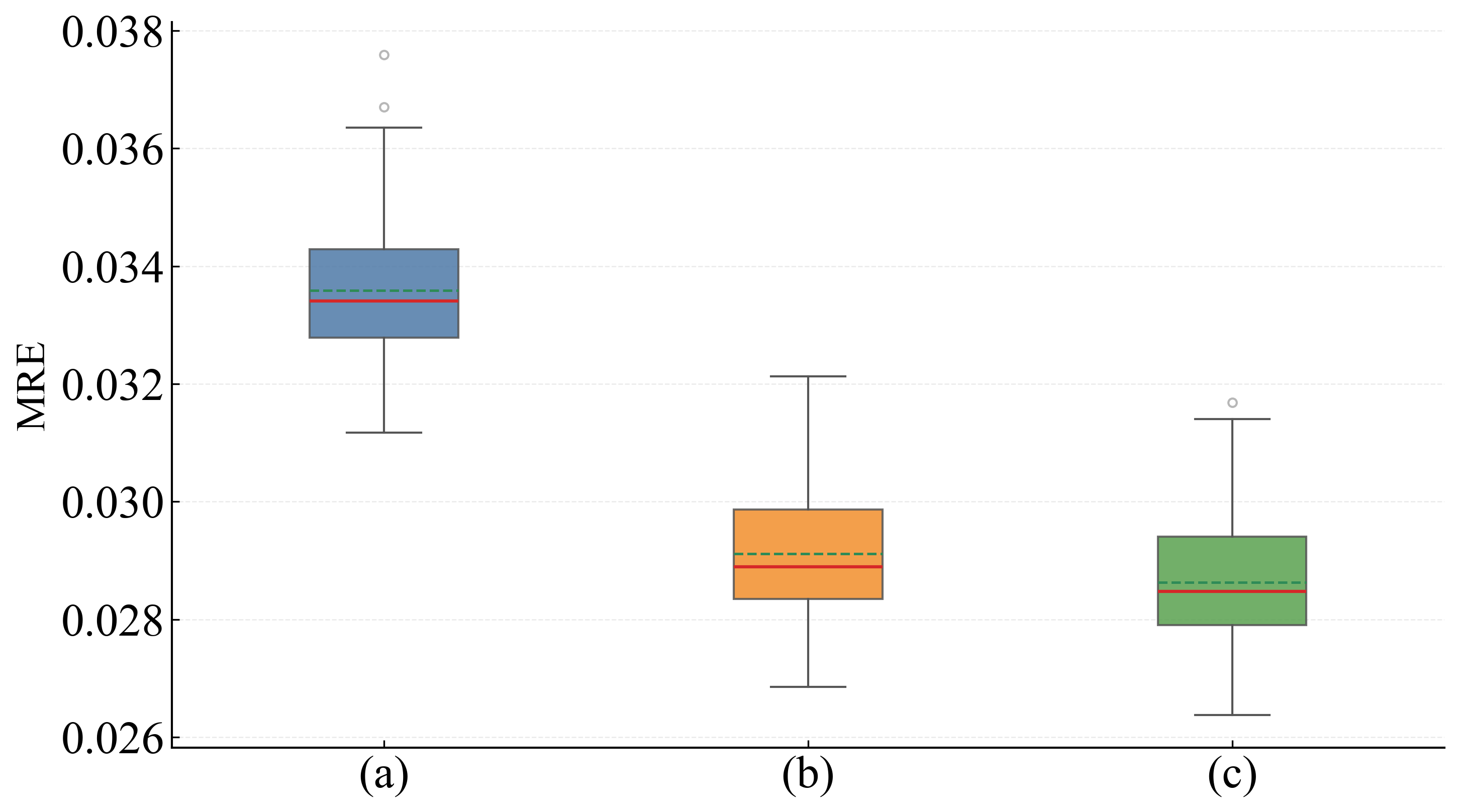}
	\caption{Prediction MRE for $41\times41\times6$ saturation fields using zero-order (a), first‑order (b) and second‑order (c) loss function. The training set consists of 200 $21\times21\times6$ and 100 $11\times11\times6$ samples.}
	\label{fig:sw41_compare_200_100}
\end{figure}

The effect of multiresolution Sobolev training on cross-resolution prediction is validated using a representative training set consists of 200 $21\times21\times6$ and 100 $11\times11\times6$ samples. The comparison using MAE averaged on both cells and samples is presented in Table~\ref{tab:FNO_Press_gd_compare of 200 and 100} for pressure and Table~\ref{tab:FNO_Sw_gd_compare of 200 and 100} for saturation. Including first-order gradients in the loss is very effective for enhancing the cross-resolution prediction accuracy, while including extra second-order gradients is still meaningful as the distribution of MRE can be improved slightly as is shown in Fig.~\ref{fig:press41_compare_200_100} for pressure and Fig.~\ref{fig:sw41_compare_200_100} for saturation.

\begin{table}[!htb]
\centering
\begin{tabular}{ccc}
\toprule
Loss Function Order &$R^2$ & MAE (MPa)\\
\midrule
0 &  0.411321 & 0.14184923 \\
1 &  0.515972 & 0.12641717 \\
2 &  0.520149 & 0.12582239 \\
\bottomrule
\end{tabular}
\caption{Prediction accuracy comparing different orders of loss function for the pressure field. The training set consists of 200 $21\times21\times6$ samples and 100 $11\times11\times6$ samples. The prediction accuracy is quantified on data of resolution $41\times41\times6$.}
\label{tab:FNO_Press_gd_compare of 200 and 100}
\end{table}

\begin{table}[!htb]
\centering
\begin{tabular}{ccc}
\toprule
Loss Function order &$R^2$ & MAE (percent)\\
\midrule
0 & 0.757315 & 1.140762 \\
1 & 0.801380 & 1.002029 \\
2 & 0.806094 & 0.987777 \\
\bottomrule
\end{tabular}
\caption{Prediction accuracy comparing different orders of loss function for the saturation field. The training set consists of 200 $21\times21\times6$ samples and 100 $11\times11\times6$ samples. The prediction accuracy is quantified on data of resolution $41\times41\times6$.}
\label{tab:FNO_Sw_gd_compare of 200 and 100}
\end{table}

\begin{figure}[!htb]
	\centering
	\vspace{0.0cm}
	\includegraphics[width=0.7\textwidth]{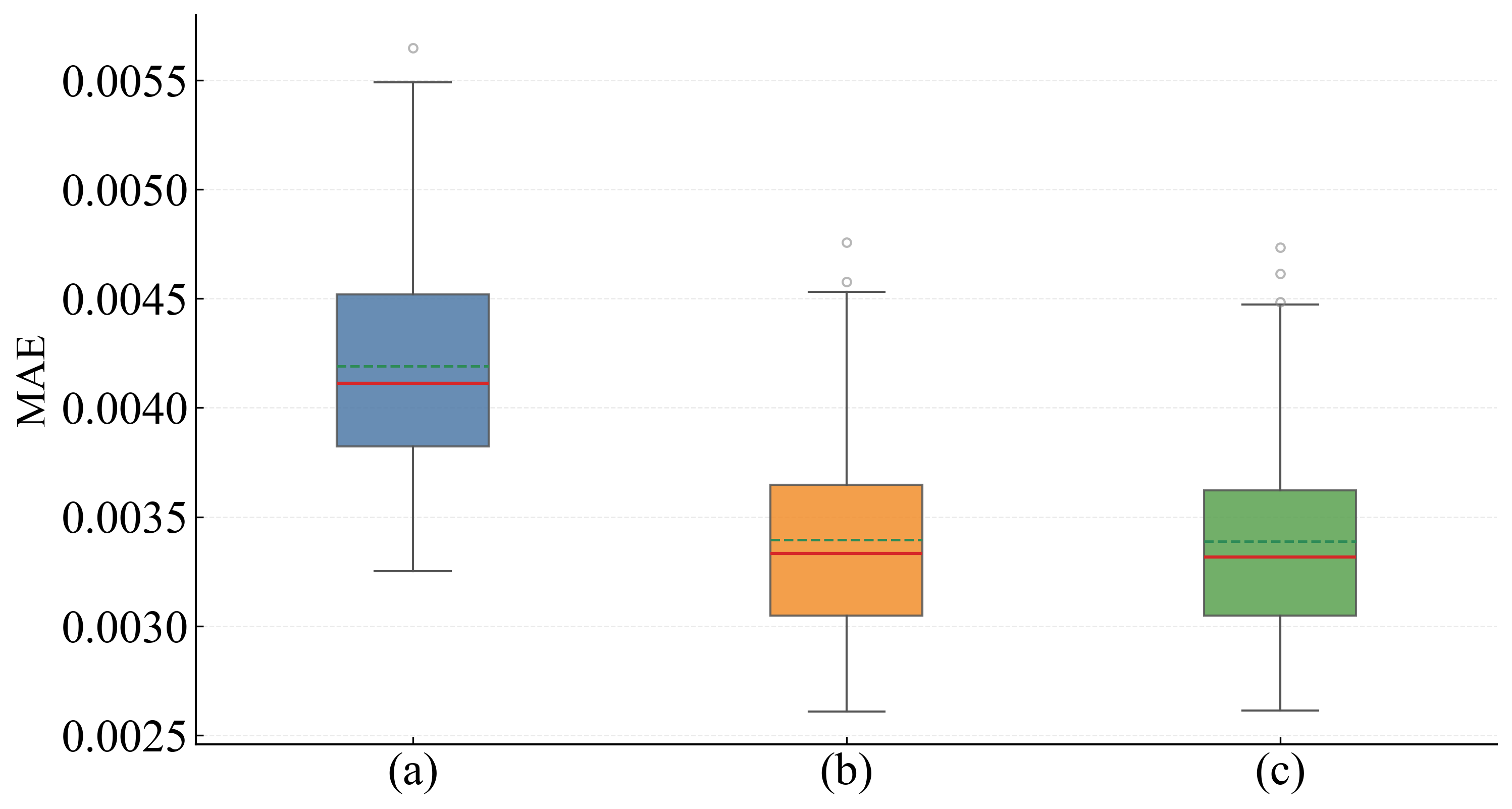}
	\caption{Prediction MRE for $41\times41\times6$ x-derivative fields of pressure using zero-order (a), first‑order (b) and second‑order (c) loss function. The training set consists of 200 $21\times21\times6$ and 100 $11\times11\times6$ samples.}
\label{fig:Press_comparison_41_MAE_boxplot_first_grad}
\end{figure}

\begin{figure}[!htb]
	\centering
	\vspace{0.0cm}
	\includegraphics[width=0.7\textwidth]{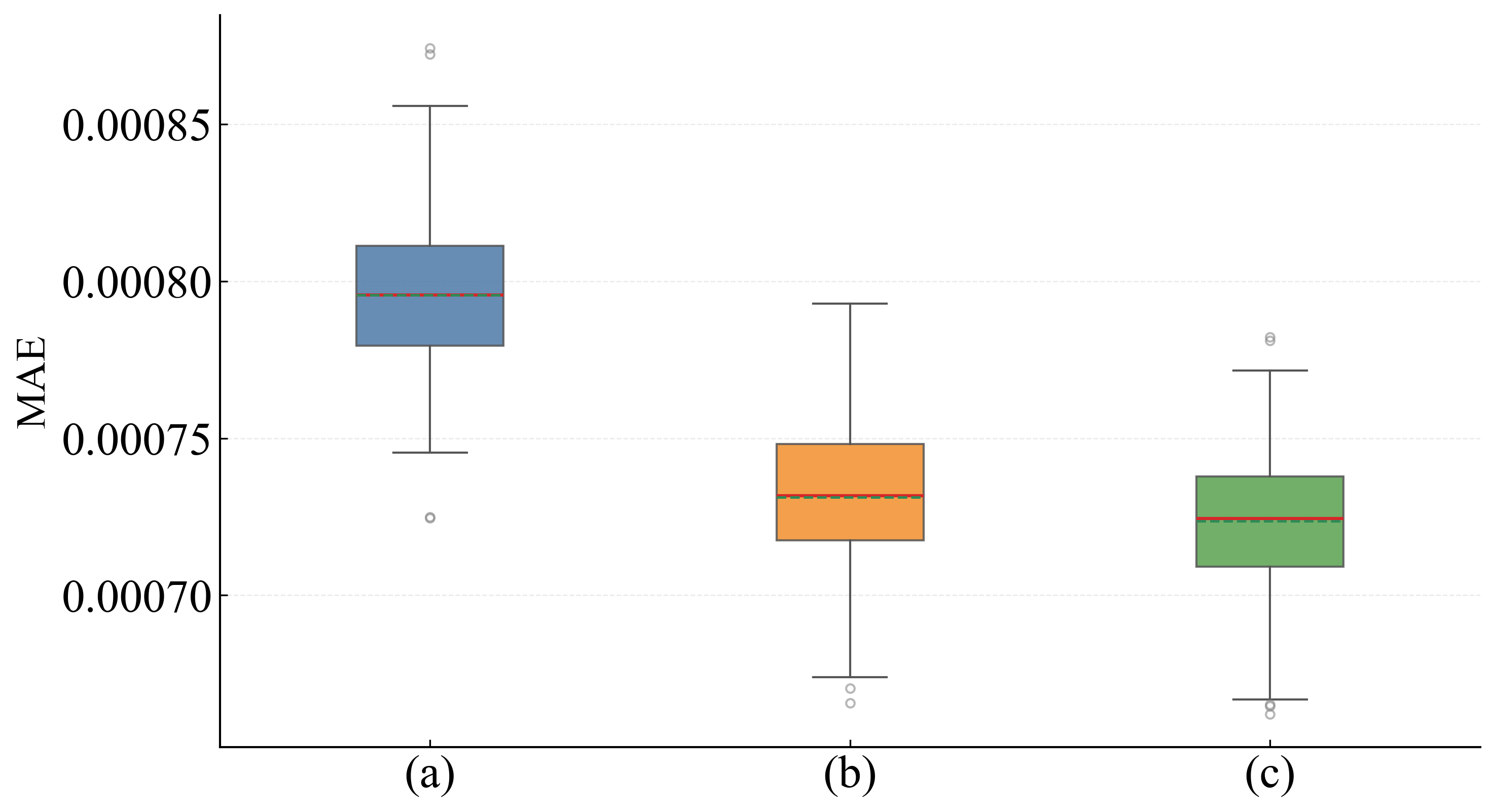}
	\caption{Prediction MRE for $41\times41\times6$ x-derivative fields of saturation using zero-order (a), first‑order (b) and second‑order (c) loss function. The training set consists of 200 $21\times21\times6$ and 100 $11\times11\times6$ samples.}
	\label{fig:Sw_comparison_41_MAE_boxplot_first_grad}
\end{figure}

Further, the effect of multiresolution Sobolev training is evaluated for enhancing the derivative fields of pressure and saturation approximated by the finite difference scheme. Fig.~\ref{fig:Press_comparison_41_MAE_boxplot_first_grad} and \ref{fig:Sw_comparison_41_MAE_boxplot_first_grad} are the box-plots showing the MRE distribution of FNO for predicting the x-derivative fields of pressure and saturation, trained without gradient terms, with first-order gradient, and with both first- and second-order gradients. Similarly as the results for the original pressure and saturation fields, including the first-order gradient terms in the loss function is very effective while including extra second-order gradient terms is still meaningful especially for saturation.

\begin{figure}[!htb]
	\centering
	\vspace{0.0cm}
	\includegraphics[width=0.99\textwidth]{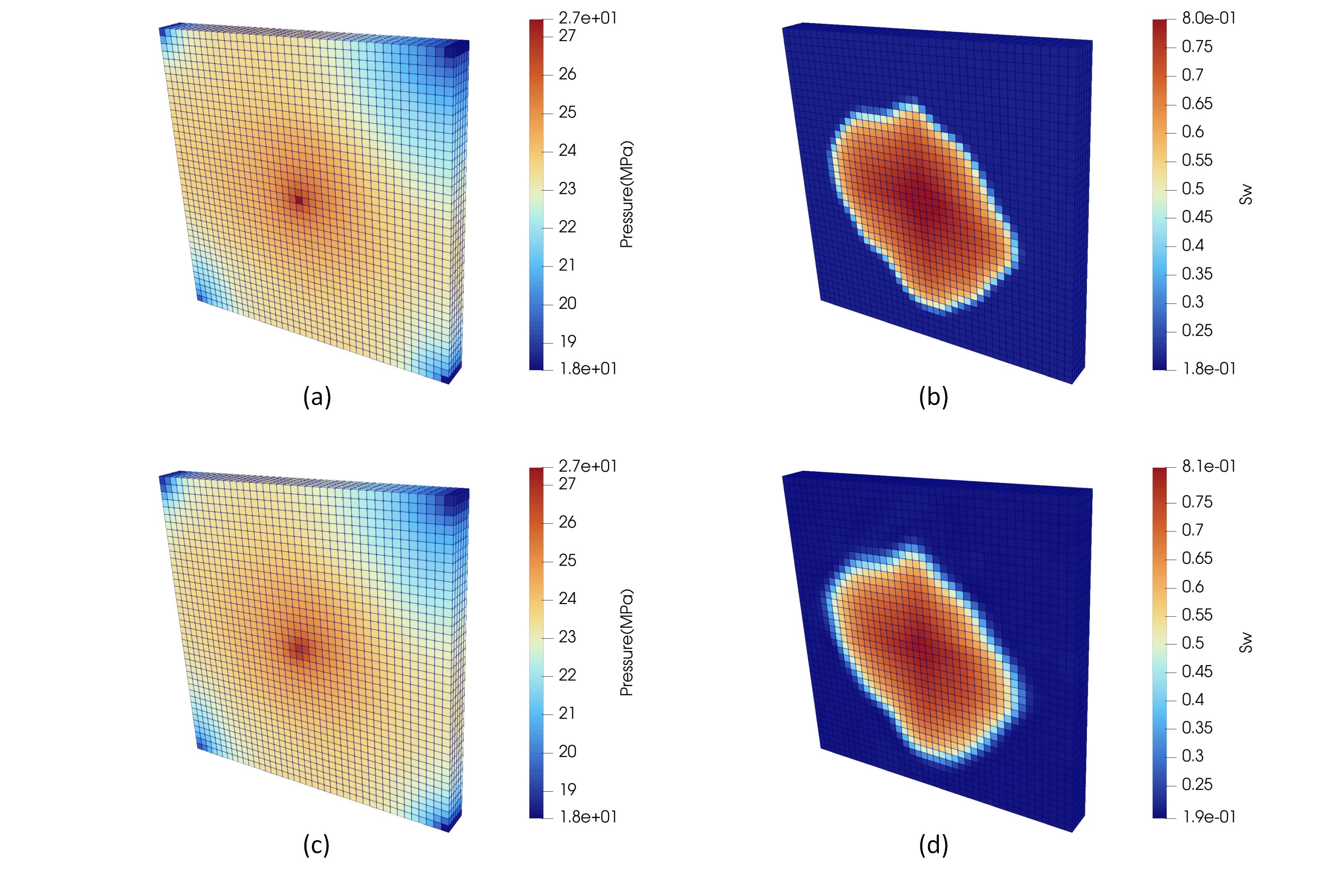}
	\caption{(a) True final time step pressure, (b) True final time step saturation, (c) Predicted final time step pressure, (d) Predicted final time step saturation. }
	\label{fig:last_time_press+sw.jpg}
\end{figure}

\begin{figure}[!htb]
	\centering
	\vspace{0.0cm}
	\includegraphics[width=0.99\textwidth]{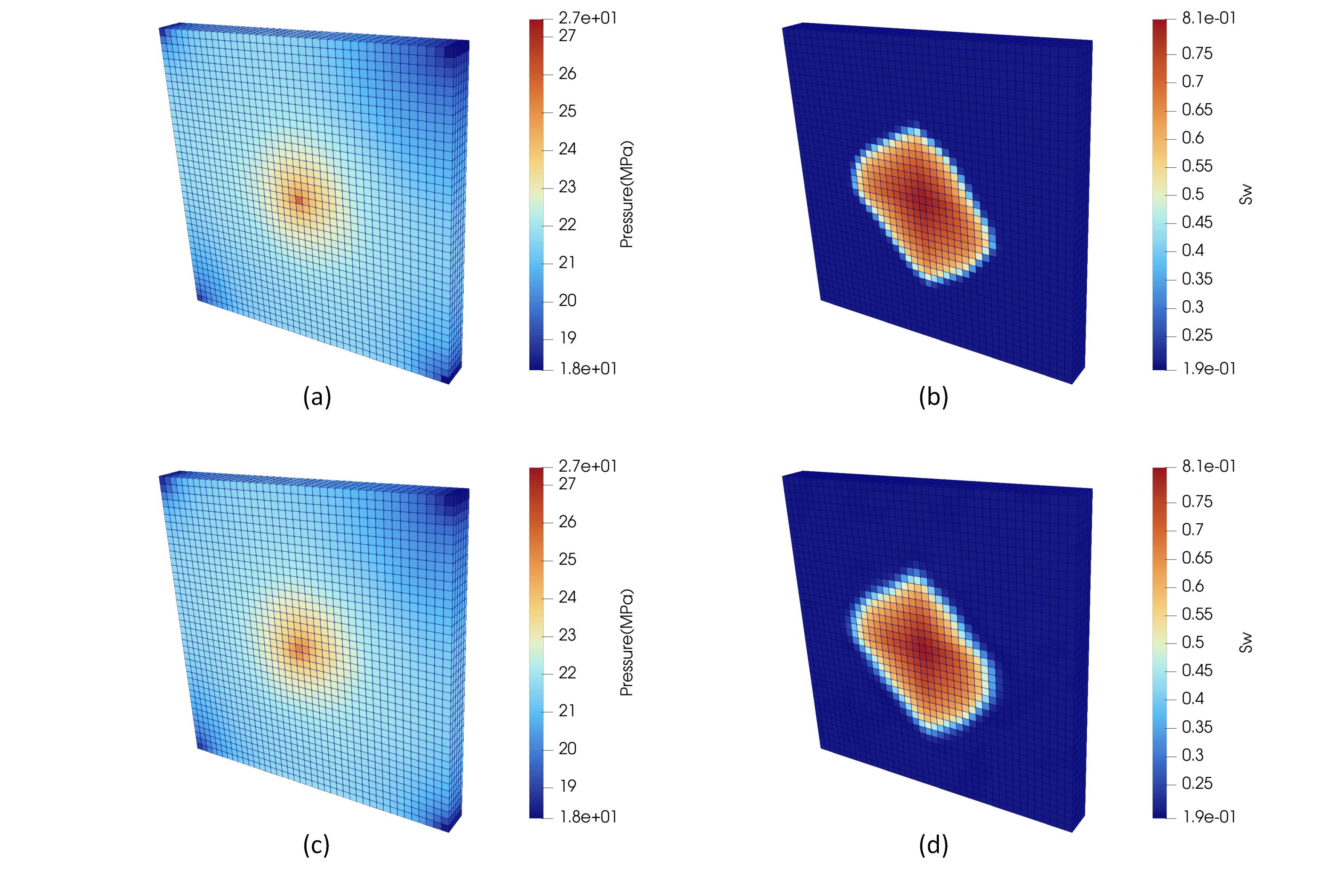}
	\caption{(a) True intermediate time step pressure, (b) True intermediate time step saturation, (c) Predicted intermediate time step pressure, (d) Predicted intermediate time step saturation. }
	\label{fig:mid_time_press+sw.jpg}
\end{figure}

The FNO-based surrogate can be used to predict cross-resolution pressure and saturation fields of different time steps.
Fig.~\ref{fig:last_time_press+sw.jpg} presents the predicted pressure and saturation fields at end of simulation at resolution level $41\times41\times6$, as compared with the ground truth for a random sample. Fig.~\ref{fig:mid_time_press+sw.jpg} presents the comparison at an intermediate time step. The predicted and true pressure fields are very close to each other visually. The average differences of swept volumes between predicted and true saturation fields are generally below 0.3\%.

\section{Conclusion}
In the current study, we have built a cross-resolution FNO-based surrogate for parametric two-phase transient compressible Darcy flows in porous media. A multiresolution Sobolev learning method is proposed to enhance the cross-resolution prediction accuracy by incorporating first and second-order gradients in the loss function. The gradients are approximated by finite difference schemes that are resolution-related and compatible with numerical simulation. The theoretical guarantee of Sobolev training for enhancing convergence have been proved in literature but with simplified assumptions that are difficult to be generalized to operator learning based on FNO. Numerical experiments have been conducted in this paper to validate the FNO-based cross-resolution surrogate and the new training scheme involving both first- and second-order gradients as well as multiresolution training data. Results show that the use of training data of different resolutions can lead to much higher training convergence as well as prediction accuracy on test samples of higher resolution. Including first-order gradients in the loss function can yield notable increase in the cross-resolution prediction accuracy while including second-order gradients can yield further but moderate increase in the cross-resolution prediction accuracy for the solution fields as well as the derivative fields. The current work could be extended to unstructured grids in future studies.

\section*{Acknowledgements}
The research is supported by the Natural Science Foundation of Shandong Province (No. ZR2024MA057), the Fundamental Research Funds for the Central Universities and the Future Plan for Young Scholars of Shandong University.

\bibliographystyle{unsrtnat}
\bibliography{references}

\end{document}